\documentclass{article}
\usepackage{amsfonts}
\usepackage{amssymb}
\usepackage{amsbsy}
\usepackage{amsmath} 
\usepackage{comment}

\newtheorem{theorem}{Theorem}
\newtheorem{lemma}[theorem]{Lemma}

\newtheorem{definition}[theorem]{Definition}

\newtheorem{corollary}[theorem]{Corollary}

\newcommand{\qed}{\hfill$\blacksquare$}
\newcommand{\halmos}{\hfill$\square$}

\newcommand{\qtwo}{\mathbf{Q_{\rm 2}}}

\newcommand{\ctwo}{\mathbf{C_{\rm 2}}}
\newcommand{\ftwo}{\mathbf{F_{\rm 2}}}

\newcommand{\zero}{\mathbf{0}}

\newcommand{\x}{z}

\newcommand{\C}{\mathbf{C}}

\newcommand{\Q}{\mathbf{Q}}

\newcommand{\Z}{\mathbf{Z}}

\newcommand{\M}{\mathcal{M}}

\newcommand{\kk}{\mathbf{k}}

\DeclareMathOperator{\col}{col}

\begin{document}

\bibliographystyle{plain}

\title{Slopes of 2-adic overconvergent modular forms with small level}
\author{L J P Kilford\footnote{Address: Mathematics 253-37, California Institute of Technology, Pasadena CA 91125.}}

\maketitle

\section{Abstract}

Let~$\tau$ be the primitive Dirichlet character of conductor~$4$, let~$\chi$ be the primitive even Dirichlet character of conductor~$8$ and let~$k$ be an integer. Then the~$U_2$ operator acting on cuspidal overconvergent modular forms of weight~$2k+1$ and character~$\tau$ has slopes in the arithmetic progression~$\left\{2,4,\ldots,2n,\ldots\right\}$, and the~$U_2$ operator acting on cuspidal overconvergent modular forms of weight~$k$ and character~$\chi \cdot \tau^k$ has slopes in the arithmetic progression~$\left\{1,2,\ldots,n,\ldots\right\}$.

We then show that the characteristic polynomials of the Hecke operators~$U_2$ and~$T_p$ acting on the space of classical cusp forms of weight~$k$ and character either~$\tau$ or~$\chi\cdot\tau^k$ split completely over~$\qtwo$.

\section{Introduction}

\begin{definition}
Let~$f$ be a cuspidal modular eigenform with $q$-expansion at~$\infty$ given by~$\sum_{n=1}^\infty a_n q^n$. Let~$f$ be \emph{normalised}; that is, $a_1=1$. The  \emph{$(p)$-slope} of~$f$ is defined to be the $p$-valuation of~$a_p$; we normalise the $p$-valuation of~$p$ to be~$1$. If we do not specify~$p$, then we mean the $2$-slope.
\end{definition}

In this paper, we prove the following theorem on the slopes of classical modular cusp forms :

\begin{theorem}
\label{introduction-slopes-theorem}
Let~$\tau$ be the nontrivial character of conductor~$4$, and let~$k$ be an integer greater than~$2$.
The slopes of the~$U_2$ operator acting on~$S_{2k-1}(\Gamma_0(4),\tau)$ are
\[
2,4,6,\ldots,2k-4.
\]
Let~$\chi$ be the even primitive Dirichlet character of conductor~$8$.
The slopes of the~$U_2$ operator acting on~$S_{k}(\Gamma_0(8),\chi \cdot \tau^k)$ are
\[
1,2,3,\ldots,k-2.
\]
\end{theorem}

As a corollary of this theorem, we also prove the following result about the field over which cusp forms of weight~$k$ and character~$\chi\cdot\tau^k$ or~$\tau$ are defined:

\begin{corollary}
\label{introduction-definition-theorem}
Let~$k$ be an integer greater than~$2$, and let~$S$ be either 
$S_{2k-1}(\Gamma_0(4),\tau)$ or $S_{k}(\Gamma_0(8),\chi \cdot\tau^k)$ 

Let~$f \in S$ be a normalised eigenform. Then the coefficients of the Fourier expansion of~$f$ are elements of~$\mathbf{Q}_2$.
\end{corollary}

This corollary gives a partial answer to an extension of Questions~4.3 and~4.4 of Buzzard~\cite{buzzard-questions}, which give a conjectural bound on the degree of the field of definition of certain spaces of modular forms over~$\qtwo$.

\section{Previous work}

Matthew Emerton determines in his thesis~\cite{emerton-thesis} the smallest slope for the 
spaces of modular cuspforms~$S_k(\Gamma_0(2^n),\chi)$, where~$\chi$ is a primitive Dirichlet character of conductor~$2^n$. 
\begin{theorem}[Emerton~\cite{emerton-thesis}, Proposition 5.1]
Let~$m$ be a positive integer greater than~$1$, and let~$\chi$ be a primitive Dirichlet character of conductor~$2^m$. The smallest slope of the $U_2$ operator acting on cuspforms of weight~$k$ and character~$\chi$ 
is~$2^{3-m}$.
\end{theorem}

If we look at the character of conductor~4 and the odd character of conductor~8, there is a CM modular form which is defined over the field~$\Q$. We quote a result of Schoeneberg, proved in Ogg~\cite{ogg}:
\begin{theorem}[Ogg~\cite{ogg}, Theorem VI.22]
Let~$i$ be the square root of~$-1$, and let~$k$ be an positive odd integer greater than~$3$ and congruent to~$1 \mod 4$. Then there is a normalised cuspidal modular eigenform in~$S_k(\Gamma_0(4),\tau)$ with $q$-expansion  
\[
f_k(q)=\frac{1}{4}\sum_{m,n \in \Z} (m+n\cdot i)^{k-1} \cdot q^{m^2+n^2}.
\]
Let~$l$ be a positive odd integer greater than~$1$. Then there is a normalised cuspidal modular eigenform in~$S_{l}(\Gamma_0(8),\tau \cdot \chi)$ with $q$-expansion 
\[
g_l(q)=\frac{1}{2}\sum_{m,n \in \Z} (m+2n\cdot i)^{l-1} \cdot q^{m^2+2n^2} 
\]
\end{theorem}
We see by inspection that the slope of~$f_k$ is~$(k-1)/2$, and that the slope of~$g_l$ is~$l-1$.

Hence we can, in certain cases, determine the smallest slope and another classical slope of~$U_2$ acting on modular newforms of level~$\Gamma_1(4)$ or~$\Gamma_1(8)$ using previously known results.

Lawren Smithline has also proved results about the slopes of classical modular forms, and the techniques used in~\cite{smithline} are similar to those in this paper.
\begin{theorem}[Smithline~\cite{smithline}, Corollary 6.1.3.3]
Let~$v$ be a non-negative integer and let~$k=2 \cdot 3^{v+1}$. Then there are exactly~$3^v$ classical modular eigenforms of weight~$k$ and level~$3$ with $3$-slope~$3^{v+1}-1$.
\end{theorem}

Buzzard and Calegari~\cite{buzzard-calegari} have proved the following theorem on the slopes of modular forms:
\begin{theorem}[Buzzard-Calegari~\cite{buzzard-calegari},~Corollary~1]
The slopes of the~$U_2$ operator on the space of cusp forms of weight~$0$ are given by
\[
1+2v_2\left(\frac{(3n)!}{n!}\right).
\]
\end{theorem}

Jacobs~\cite{dan-jacobs-thesis} has proved, using similar techniques to those in this paper, that the slopes of the operator~$U_3$ acting on the space of automorphic forms over a definite quarternion algebra are in the arithmetic progression~$\left\{1/2,3/2,5/2,\ldots\right\}$.

\section{Overconvergent modular forms}

A famous quote of Jacques Hadamard~\cite{hadamard} says that ``the shortest and best way between two truths of the real domain often passes through the imaginary one." It seems that often the best way to prove results like Theorem~\ref{introduction-slopes-theorem} about classical modular forms is to prove a theorem for the overconvergent modular forms and then derive the theorem for classical modular forms as a consequence.

We therefore recall the definition of the 2-adic overconvergent modular forms, first by defining overconvergent modular forms of weight~$0$, and then by deriving the definition for forms with weight and character.

Following Katz~\cite{katz}, section~2.1, we recall that, for~$E$ be an elliptic curve over an $\mathbf{F}_2$-algebra~$R$, there is a mod~$2$ modular form~$A(E)$ called the \emph{Hasse invariant}, which has the $q$-expansion over~$\mathbf{F}_2$ equal to~$1$.

We consider the Eisenstein series of weight~$4$ and tame level~1 defined over~$\mathbf{Z}$, with~$q$-expansion
\[
E_4(q):=1+240 \sum_{n = 1}^\infty \left( \sum_{0 < d | n} d^3 \right) \cdot q^n.
\]

We see that~$E_4$ is a lifting of~$A(E)^4$ to characteristic~0, as the reduction of~$E_4$ to characteristic~$2$ has the same $q$-expansion as~$A(E)^4$, and therefore~$E_4 \mod 2$ and~$A(E)^4$ are both modular forms of level~1 and weight~4 defined over~$\mathbf{F}_2$, with the same $q$-expansion.

Now the value of~$E_4(E)$ may not be well-defined, but it can be shown that the valuation~$v_2(E_4(E))$ is well-defined. This will allow us to define the ordinary locus of~$X_0(2^m)$ and certain neighbourhoods of it. We follow the work of Coleman~\cite{coleman}, and first define structures on the modular curve~$X_1(2^m)$.

\begin{definition}[Coleman~\cite{coleman}, page~448]
Consider~$X_1(2^m)_{/\mathbf{Q}_2}$ as a rigid analytic space, and let~$t$ be a point of~$X_1(2^m)$.

If~$t$ is a point of~$X_1(2^m)$ which corresponds to a cusp, then we define~$v(E_4(t))=0$, following~\cite{buzzard-analytic}, section~4.

We define the ordinary locus of~$X_1(2^m)$ to be the set of points~$t$ of~$X_1(2^m)$ such that~$v_2(E_4(t))=0$, and define~$Z_1(2^m)$ to be the rigid connected component of the ordinary locus in~$X_1(2^m)$ which contains the cusp~$\infty$. This is a rigid analytic space.


\end{definition}
In~\cite{coleman-mazur}, page~36, it is shown that~$Z_1(2^m)$ is an affinoid subdomain of the rigid space~$X_1(2^m)_{/\mathbf{Q}_2}$.

We will perform calculations in later sections on the modular curve~$X_0(2^m)$, which we will now define.
\begin{definition}
Consider~$X_1(2^m)$ as a modular curve. We see that the group~$G:=(\mathbf{Z}/2^m \mathbf{Z})^\times$ acts upon the non-cuspidal points of~$X_1(2^m)$, by the following action: if~$a \in (\mathbf{Z}/2^m \mathbf{Z})^\times$, then the action of~$a$ sends the pair~$(E,P)$ to~$(E,aP)$. This action extends to the cuspidal points of~$X_1(2^m)$, and it sends cusps to cusps.

 We will define the modular curve~$X_0(2^m)_{/\mathbf{Q}_2}$ to be the quotient of~$X_1(2^m)$ by~$(\mathbf{Z}/2^m \mathbf{Z})^\times$.
\end{definition}

We note that the action of the group~$G$ does not change the valuation of a given elliptic curve~$E$. We define~$Z_0(2^m)$ to be the rigid connected component of the ordinary locus in~$X_0(2^m)$ which contains the cusp~$\infty$. It is a rigid analytic space.

We will now define strict affinoid neighbourhoods of~$Z_0(2^m)$.
\begin{definition}[Coleman~\cite{coleman}, Section~B2]
\label{connected-component}
We think of~$X_0(2^m)$ as a rigid space over~$\mathbf{Q}_2$, and we let~$t \in X_0(2^m)(\overline{\mathbf{Q}}_2)$ be a point, corresponding either to an elliptic curve defined over a finite extension of~$\Q_2$, or to a cusp.
Let~$w$ be a rational number, such that~$0 < w < \min(2^{2-m}/3,1/4$).

We define~$Z_0(2^m)(w)$ to be the connected component 
of the affinoid 
\[
\left\{t \in X_0(2^m): \; v_2(E_4(t)) \le 4w\right\}
\]
which contains the cusp~$\infty$. 
\end{definition}
The condition involves~$4w$ rather than~$w$ because we are working with a lifting of the fourth power of the Hasse invariant. Also, note that as~$E_4$ is a lifting of the mod~$2$ modular form~$A^4$, and that any another lifting of~$A^4$ would be of the form~$E_4+2F$, where~$F$ is a modular form, then this valuation is well-defined if~$0\le v(E_4(t)) < 1$. This corresponds to the condition~$0 \le w < 1/4$ in Definition~\ref{connected-component}.

We can now define overconvergent modular forms of weight~$0$.

\begin{definition}[Coleman,~\cite{coleman-overconvergent}, page~397]
Let~$w$ be a rational number, such that~$0 < w < \min(2^{2-m}/3,1/4)$. Let~$\mathcal{O}$ be the structure sheaf of~$Z_0(2^m)(w)$.

We call sections of~$\mathcal{O}$ on~$Z_0(2^m)(w)$ 
\emph{$w$-overconvergent 2-adic modular forms of weight~$0$ and level~$\Gamma_0(2^m)$}.

If a section~$f$ of~$\mathcal{O}$ is a $w$-overconvergent modular form, then we say that~$f$ is an \emph{overconvergent $2$-adic modular form}.

Let~$K$ be a complete subfield of~$\mathbf{C}_2$, and define~$Z_0(2^m)(w)_{/K}$ to be the affinoid over~$K$ induced from~$Z_0(2^m)(w)$ by base change from~$\qtwo$. The space~$M_0(2^m,w;K):=\mathcal{O} (Z_0(2^m)(w)_{/K})$ of $w$-overconvergent modular forms of weight~$0$ and level~$\Gamma_0(2^m)$ is a $K$-Banach space.
\end{definition}

We now use non-cuspidal modular forms of the desired weight and character to define overconvergent modular forms with non-zero weight.

\begin{definition}
Let~$w$ be a real number such that~$0 < w < \min(2^{2-m}/3,1/4)$. Let~$k$ be an integer and let~$\chi$ be a character such that~$\chi(-1)=(-1)^k$, and let~$K$ be a complete subfield of~$\ctwo$. Let~$E_{k,\chi}$ be the Eisenstein series of weight~$k$ and character~$\chi$.


The space of overconvergent 2-adic modular forms of weight~$k$ and character~$\chi$ is given by
\[
\M_{k,\chi}(2^m,w;K):=E^*_{k,\chi} \cdot \M_0(2^m,w;K).
\]
This is a Banach space over~$K$.
\end{definition}

There are continuous Hecke operators~$T_p$ which act on~$\M_{k,\chi}(2^m,w;K)$; the operator~$U:=U_2$, which is defined on $q$-expansions as 
\[
U\left(\sum_{n=0}^\infty a_n q^n\right) = \sum_{n=0}^\infty a_{2n} q^n
\]
is also \emph{compact} and therefore has a spectral theory.

As a consequence of results of Coleman, we have the following theorem:
\begin{theorem}[Coleman~\cite{coleman}, Theorem~B3.2]
Let~$w$ be a real number such that~$0 < w < \min(2^{2-m}/3,1/4)$, let~$k$ be an integer and let~$\chi$ be a character such that~$\chi(-1)=(-1)^k$.

The characteristic polynomial of~$U_2$ acting on overconvergent 2-adic modular forms of weight~$k$ and character~$\chi$ is independent of the choice of~$w$.
\end{theorem}
This theorem allows us to choose a convenient value of~$w$ and prove results for that~$w$, and guarantees that these results will hold for any~$w$.

The connected component in Definition~\ref{connected-component} is hard to work with. We will therefore rewrite it in terms of modular functions of level greater than~$1$, to prove the following theorem:

\begin{theorem}
The spaces of overconvergent modular forms of weight~0 and level~$N=4 \text{ or }8$ are Tate algebras in one variable over~$\qtwo(2^{4/N})$.
\end{theorem}

We have given a valuation on the points~$t$ of the rigid space~$X_0(2^m)$, based on the lifting of the Eisenstein series~$E_4$. We recall that the modular $j$-invariant is defined to be~$j=E_4^3/\Delta$. Therefore, we see that, if the elliptic curve corresponding to~$t$ has good reduction, then~$\Delta(t)$ has valuation~0, and therefore that 
\[
v(t)=\frac{1}{4} v(E_4(t))=\frac{1}{12} v((E_4)^3(t))=\frac{1}{12} v(j(t)).
\]
From Lemma~2.3 of~\cite{emerton-thesis}, we see that there is a modular function~$j_8$ which is a uniformiser on~$X_0(8)$. It has $q$-expansion at~$\infty$
\[
j_8=\frac{1}{q\prod_{n=1}^\infty (1+q^n)^4(1+q^{2n})^2(1+q^{4n})^4}=\left(\frac{\Delta(q) \Delta(q^4)}{\Delta(q^2) \Delta(q^8)}\right)^{1/12}.
\]
Also, $j_8(\infty)=\infty$.

There is another modular function~$j_{16}$ which is a uniformiser on~$X_0(16)$, with $q$-expansion at~$\infty$ given by
\[
j_{16}=\frac{1}{q \prod_{n=1}^\infty (1+q^n)^2 (1+q^{2n})(1+q^{4n})(1+q^{8n})^2}=\left(\frac{\Delta(q^{16})^2 \Delta(q^2)}{\Delta(q^8)\Delta(q)^2}\right)^{1/24}.
\]
We see also that~$j_{16}(\infty)=\infty$.

By an explicit calculation of $q$-expansions, using the formulae in Chapter~2 of~\cite{emerton-thesis}, 
we see that
\[ 
j=\frac{(j_8^4+256j_8^3+5120j_8^2+32768j_8+65536)^3}{(j_8^2+16j_8+64)\cdot(j_8+4)}
\]
and
\[
\frac{1}{j_8}=\frac{1}{j_{16}}+\frac{2}{j^2_{16}}.
\]

Because we know that~$j_8(\infty)=\infty$, the connected component of~$Z_0(8)$ which contains~$\infty$ is of the form~$v(j_8) < D$ for some rational number~$D$. We see that, if~$v(j_8) < 2$, then~$v(j_8)=v(j)$. This means that we have shown that
\[
Z_0(8)(w)=\left\{x \in X_0(8):\;v(j_8(x)) \le 12w\right\} \text{ for }0 < w < 1/4.
\]
Similarly, we see that the connected component of~$Z_0(16)$ which contains~$\infty$ is of the form~$v(j_8) < D$ for some rational number~$D$. We see that, if~$v(j_{16}) < 1$, then~$v(j_8)=v(j_{16})$, and therefore~$v(j_{16})=v(j_8)$. This means that we have shown that
\[
Z_0(16)(w)=\left\{x \in X_0(8):\;v(j_{16}(x)) \le 12w\right\}, \text{ for }0 < w < 1/6.
\]

We now define another modular function on~$X_0(2^m)$, in terms of Eisenstein series.

Let~$k$ be an integer, and let~$\theta:(\Z/n\Z)^\times \rightarrow \C^\times$ be a primitive Dirichlet character, such that~$\theta(-1)=(-1)^k$. Recall from Washington~\cite{washington}, page~30, that the \emph{extended Bernoulli numbers}~$B_{k,\theta}$ are defined by
\begin{displaymath}
\sum^{N}_{a=1} \frac{\theta(a)\cdot t \cdot \exp({at})}{\exp({Nt})-1} = \sum_{k=0}^\infty B_{k,\theta} \cdot 
\frac{t^i}{i!} 
\end{displaymath}

We define the Eisenstein series~$E^*_{k,\theta}$ to be
\[
E^*_{k,\theta}:= \frac{-B_{k,\theta}}{2k}+ \sum^{\infty}_{n = 1} \left( \underset{(d,p)=1}{\sum_{0 < d | n}} \theta(d) \cdot d^{k-1} \right) \cdot q^n,
\]
where~$B_{k,\theta}$ is the extended Bernoulli number attached to~$k$ and~$\theta$.

There is an operator~$V$ on the space of modular forms~$S_k(\Gamma_0(N))$; its effect on $q$-expansions is to send~$q$ to~$q^2$. We define~$V^*_k$ to be~$V(E^*_k)$.

We define modular functions on~$X_0(8)$ and~$X_0(16)$ by
\[
z_4:=\frac{E^*_{1,\tau}/V^*_{1,\tau}-1}{2} = \frac{2}{j_8+2}
\]
and
\[
z_8:=\frac{E^*_{1,\chi\tau}/V^*_{1,\chi\tau}-1}{\sqrt{2}} = \frac{\sqrt{2}}{j_{16}+2},
\]
where we choose and fix a square root of~$2$ in~$\ctwo$.
These identities can be verified by explicit calculation.

Let~$w$ be a rational number such that~$0 < w < 1/6$. Then using the formulae above, we see that
\[
Z_0(8)(w)=\left\{x \in X_0(8):\;v(z_4(x)) \ge 1-12w\right\}.
\]
We now choose~$w=1/12$, to obtain
\[
Z_0(8)(1/12)=\left\{x \in X_0(8):\;v(z_4(x)) \ge 0\right\}.
\]
Now, the rigid functions on the closed disc over~$\mathbf{Q}_2$ with centre~0 and radius~1 are defined to be power series of the form
\[
\sum_{n \in \mathbf{N}}a_n z^n\;:\;a_n \in \mathbf{Q}_2,\;a_n \rightarrow 0.
\]
Therefore, the $1/12$-overconvergent modular forms of level~$\Gamma_0(4)$ and weight~0 are 
\[
\mathbf{Q}_2\langle z_4 \rangle,
\]
which is what we wanted to show.
We now follow the same procedure for~$X_0(16)$.
Let~$w$ be a rational number such that~$0 < w < 1/12$. Then using the formulae above, we see that
\[
Z_0(16)(w)=\left\{x \in X_0(16):\;v(z_8(x)) \ge 1/2-12w\right\}.
\]
We now choose~$w=1/24$, to obtain
\[
Z_0(16)(1/24)=\left\{x \in X_0(16):\;v(z_8(x)) \ge 0\right\}.
\]
Now, the rigid functions on the closed disc over~$\mathbf{Q}_2$ with centre~0 and radius~1 are defined to be power series of the form
\[
\sum_{n \in \mathbf{N}}a_n z^n\;:\;a_n \in \mathbf{Q}_2(\sqrt{2}),\;a_n \rightarrow 0.
\]
Therefore, the $1/24$-overconvergent modular forms of level~$\Gamma_0(8)$ and weight~0 are 
\[
\mathbf{Q}_2(\sqrt{2})\langle z_8 \rangle,
\]
so we have shown that these spaces of modular forms are Tate algebras in one variable.\qed




We now show that the odd powers of~$z_N$ are sent to~$0$ under the~$U_2$ operator.
\begin{lemma}
Let~$N=4$ or~$8$ and let~$i$ be a positive integer. Then
\[
U_2(z_N^{2i+1})=0,\text{ and } U(z_N^{2i})=\left(U(z_N^2)\right)^{i}.
\]
\end{lemma}

Recall that the Eisenstein series~$E^*_\kk$ is an eigenvector with eigenvalue~1 for~$U_2$, and that $U_2(V(E^*_\kk))=E^*_\kk$. Then we see that we have (for~$\mu=2 \text{ or } \sqrt{2}$):
\begin{eqnarray*}
U_2(\x)&=&U_2\!\!\left(\frac{E^*_\kk/V^*_\kk-1}{\mu}\right)=\frac{1}{\mu} \cdot U_2\!\!\left(\frac{E^*_\kk-V^*_\kk}{V^*_\kk}\right) \\ &=&\frac{1}{\mu E^*_\kk}\cdot U_2(E^*_\kk-V^*_\kk)=\frac{1}{\mu E^*_\kk}\cdot (E^*_\kk-E^*_\kk) =0.
\end{eqnarray*}
Hence we see that~$z_N$ has only odd $q$-coefficients, and that therefore~$z_N^2$ has only even $q$-coefficients.

Therefore~$z_N^{2i+1}$ has only odd $q$-coefficients. Hence for all non-negative integers~$t$, we see that
\[
U_2\skew{-2}{\left(z_N^{2t+1}\right)=0}.
\]

Because we have just shown that~$z_N$ has only odd $q$-coefficients, we see that
\[
z_N=q F(q^2)=q V(F(q)),
\]
for some power series~$F(q)$.

Therefore we have
\[
U_2(z_N^{2i})=U_2(q^{2i} V(F(q)^{2i}))=U_2(V(q^i F(q)^{2i})),
\]
and hence we see that
\[
U_2(z_N^{2i})=q^i F(q)^{2i}=\left(q F(q)^2 \right)^i=U_2(z_N^2)^i,
\]
which proves the Lemma. \qed

Because we have written down the overconvergent modular forms as an explicit Banach space, we can write down its \emph{Banach basis}: the set~$\left\{z_4,z_4^2,z_4^3,\ldots\right\}$ forms a Banach basis for the overconvergent modular forms of weight~$0$ and level~$\Gamma_0(4)$ and the functions~$\left\{z_8,z_8^2,z_8^3,\ldots\right\}$ form a Banach basis for the overconvergent modular forms of weight~$0$ and level~$\Gamma_0(8)$.

These Banach bases are composed of weight~$0$ modular functions --- we want to be able to consider the action of the~$U$ operator on overconvergent modular forms with non-zero weight-character~$\kk$. Using an observation from the work of Coleman, we will be able to move between weight-character~$\zero$ and weight-character~$\kk$ via multiplication by a suitable quotient of modular forms.

From the discussion in Coleman~\cite{coleman}, page~450, we see that the~$U_2$ operator acting on overconvergent modular forms of weight-character~$(k\cdot t,\theta^t)$ where~$t$ is odd has the same characteristic power series as the composition of the $U$ operator acting on overconvergent modular forms of weight-character~$\zero$ with multiplication by the~$z$-expansion of~$(E^*_{k,\theta}/V^*_{k,\theta})^t$.

For level~$4$ and weight~$2t+1$ and level~$8$ and weight either of the form~$2t+1$ or~$4k+2$, we can write the weight~$k$ multiplier as a power of~$1+2z_4$ or~$1+\sqrt{2}z_8$. If~$N=8$ and the weight~$k$ divides~$4$, then we can write~$E^*_{k,\theta}/V^*_{k,\theta}$ as a rational function of~$z_8$, because~$z_8$ is a rational function of the uniformiser~$j_{16}$ of the genus~$0$ modular curve~$X_0(16)$. For instance, we see that
\[
\frac{E^*_{4,\chi}}{V^*_{4,\chi}}=\frac{11+2z_8+24z_8^2-48z_8^3-16z_8^4-352z_8^5}{11+24z_8^2-16z_8^4}.
\]

We can consider the action of the~$U_2$ operator on these spaces of overconvergent modular forms.

\begin{definition}
Let~$M=(m_{i,j})$ be the infinite compact matrix of the operator~$U_2 \circ (E^*_\kk/V^*_\kk)^t$ acting on the Banach basis~$\left\{z_N,z_N^2,\ldots\right\}$, where~$m_{i,j}$ is defined to be the coefficient of~$z_N^j$ in the $z_N$-expansion of~$U_2(z_N^i) \cdot (E^*_\kk/V^*_\kk)^t$.
\end{definition}

Because we know that~$\M_0$ is a Banach space, we can show that the matrix~$M$ is a compact matrix; in other words, the trace, determinant and characteristic power series of~$M$ are all well-defined.

We will use a theorem of Serre to prove our theorem on the slopes of~$U$ acting on~$\M_\kk$.
\begin{theorem}[Serre~\cite{serre}, Proposition~7]
\label{serre-proposition}

\begin{enumerate}
\item Let~$M$ be an~$n \times n$ matrix defined over a finite extension of~$\Q_2$, and let~$0 \ne r \in \Q$. Let~$\mathop{det}(1-tM)=\sum_{i= 0}^n c_i t^i$. Let~$M_m$ be the matrix formed by the first~$m$ rows and columns of~$M$.

Assume that
\begin{enumerate}
\item For all positive integers~$m$ such that~$1 \le m \le n$, the valuation of~$\mathop{det}(M_m)$ is~$r \cdot \frac{m(m+1)}{2}$.

\item The valuation of elements in column~$j$ is at least~$r \cdot j$.
\end{enumerate}

Then we have that, for all positive integers~$m$ such that~$1 \le m \le n$, $v_2(c_m)=r \cdot \frac{m(m+1)}{2}$.

\item Let~$M_\infty$ be a compact infinite matrix. If there is a sequence of finite matrices~$M_m$ which tend to ~$M_\infty$, then the finite characteristic power series~$\mathop{det}(1-tM_m)$ converge coefficientwise to~$\mathop{det}(1-tM_\infty)$, as~$m \rightarrow \infty$.

\end{enumerate}
\end{theorem}


We now quote a result of Coleman that tells us that overconvergent modular forms of small slope are in fact classical modular forms:

\begin{theorem}[Coleman~\cite{coleman-overconvergent}, Theorem~1.1]
\label{coleman-ocgt}
Let~$k$ be a non-negative integer. Every $2$-adic overconvergent modular eigenform of weight~$k$ 
with slope strictly less than~$k-1$ is a classical modular form.
\end{theorem}

We will now state the theorem on the slopes of overconvergent modular forms of weight-character~$(2k+1,\tau)$ and weight-character~$(k,\chi \cdot \tau^k)$. This, combined with Theorem~\ref{coleman-ocgt} and a knowledge of the dimensions of spaces of classical cusp forms, will suffice to prove Theorem~\ref{introduction-slopes-theorem} and Corollary~\ref{introduction-definition-theorem}.
\begin{theorem}
\label{main-theorem}
Let~$k$ be an integer, let~$\tau$ be the primitive Dirichlet character of conductor~$4$ and let~$\chi$ be the even primitive Dirichlet character of conductor~$8$.
\begin{enumerate}

\item The slopes of overconvergent 2-adic cuspidal eigenforms of weight~$2k+1$ and character~$\tau$ are~$\left\{2 i \right\}_{i \in \mathbf{N}}$.

\item The slopes of overconvergent 2-adic cuspidal eigenforms of weight~$k$ and character~$\chi \cdot \tau^k$ are~$\left\{i \right\}_{i \in \mathbf{N}}$.

\end{enumerate}
\end{theorem}

We now recall a theorem of Cohen and Oesterl\'e:
\begin{theorem}[Cohen-Oesterl\'e~\cite{cohen-oesterle}, Th\'eor\`eme~1]
\label{dimensions-of-mf}
Let~$\chi$ be a primitive Dirichlet character of conductor~$2^m > 2$, and let~$k$ be a positive integer.
Let~$\kk=(k,\chi)$ be an integral weight-character.

The dimension of the space of cuspidal modular forms of weight-character~$\kk$ is
\[
2^{m-3} (k-1)-1.
\]
\end{theorem}

We see that the slopes of the first~$2^{m-3} (k-1)-1$ overconvergent modular forms of level~$\Gamma_0(2^m)$ and primitive Dirichlet character are
\[
2^{3-m},2^{2(3-m)},\ldots,k-1-2^{3-m}.
\]
Therefore, using Theorem~\ref{coleman-ocgt}, we see that all of these slopes are classical, because~$k-1-2^{3-m} < k-1$. Hence we have proved Corollary~\ref{introduction-definition-theorem}.
\qed
\section{Necessary conditions on matrices}

We will show that we can apply Theorem~\ref{serre-proposition} by proving the following theorem: 
\begin{theorem}
\label{determinant-theorem}
Let~$N=4 \text{ or }8$, and let~$M$ be the matrix of the $U$-operator acting on the Banach basis~$\left\{(z_N)^i\right\}$. Define the set~$a_i$ by $U(z_N^2):=\sum_{i=1}^\infty a_i (z_N)^i$. Assume that
\[
\diamondsuit:\;v(a_i)=i \cdot 4/N,\; \mathrm{for\;} i\;\mathrm{odd},\;\mathrm{and}\; v(a_i) > i \cdot 4/N,\;{\rm for\;} i\; {\rm even}.
\]
Then the valuation of the determinant of the~$i \times i$ matrix~$M_n$ is~$8i/N$.
\end{theorem}

We can show the precondition of this theorem directly, by considering the two identities of modular functions
\[
U(z_4^2)=\frac{2z_4}{(1+2z_4)^2}\text{ and }U(z_8^2)=\frac{\sqrt{2}z_8}{1+2z_8^2},
\]
where we have chosen a square root of~$2$ in the extension~$\qtwo(\sqrt{2})$. These can be verified by multiplying out both sides and transforming them into an identity of modular \emph{forms}; we then use a theorem of Sturm~\cite{sturm} to show that both sides of the equation are the same by checking the $q$-expansions at~$\infty$.

We then view these rational functions of~$z_N$ as generating functions for power series in~$z_N$.

We will also prove a result about the ring over which the $q$-expansions at~$\infty$ of the cusp forms of weight-character~$\kk$ are defined. This will allow us to prove Corollary~\ref{introduction-definition-theorem}.
\begin{corollary}
\label{field-of-definition}
Let~$\kk=(k,\chi)$ be an integral weight-character, where~$\chi$ is a primitive Dirichlet character of conductor~$4$ or~$8$ and~$k$ is a positive integer such that~$\chi(-1)=(-1)^k$. Let~$K$ be the field~$\mathbf{Q}_2$ if~$N=4$ or~$\mathbf{Q}_2(\sqrt{2})$ if~$N=8$.

Let~$f=\sum_{n=1}^\infty a_n q^n$ be a normalised classical cuspidal modular eigenform of weight-character~$\kk$.

Then $a_n \in K$ for all~$n$.
\end{corollary}

Theorem~\ref{main-theorem} tells us that the slopes of overconvergent modular forms of weight-character~$\kk$ are distinct. Hence by applying Theorem~\ref{dimensions-of-mf} and Theorem~\ref{coleman-ocgt} we see that the slopes of classical modular forms of weight-character~$\kk$ are also distinct.

We will now recall a fact from Ribet~\cite{ribet-nebentypus}, page~21, to prove this Corollary. Let~$\sigma$ be an element of~$\mathop{Gal}(\overline{K}/K)$. Then we have that
\[
\sigma(f):=\sum_{n=1}^\infty \sigma(a_n) q^n
\]
is a classical cuspidal modular eigenform of weight-character~$\kk$. 

We see that the valuation of~$\sigma(a_2)$ is the same as that of~$a_2$, because the characteristic polynomial of~$a_2$ is stable under conjugation by~$\sigma$. Therefore, $\sigma(f)$ is an eigenform of weight-character~$\kk$ with the same slope as~$f$. Hence~$\sigma(f)=f$, because there is only one classical eigenform of weight-character~$\kk$ which has any given slope, by Theorem~\ref{coleman-ocgt} and Theorem~\ref{main-theorem}.

This means that~$\sigma(f)=f$ for all~$\sigma$. Therefore~$a_n \in K$ for all positive integers~$n$.\halmos

We will perform a series of transformations to obtain a matrix which we can apply Theorem~\ref{serre-proposition} to.

We first note that the odd-numbered columns of the matrix~$M_n$ are identically zero, because we have shown that~$U(z_N^{2a+1})=0$. 

We consider the matrix~$O_n$, defined by
\[
\left(O_n\right)_{i,j}:=\left(M^\prime_{2n}\right)_{2i,2j},\;{\rm where}\; 1 \le i,j \le n.
\]
This has the same characteristic power series as~$M_{2n}$.

We now show that the matrices~$O_n$ have determinant of valuation~$8\cdot n/N$, in order to be able to use Theorem~\ref{serre-proposition}.

To do this we will pre- and post-multiply the matrix~$O_n$ by diagonal matrices to obtain a matrix~$O^\prime$ which has elements of valuation at least~$8i/N$ in column~$i$. Let~$D(\alpha)$ be the diagonal matrix with $(i,i)^{th}$ coefficient~$\alpha^i$. We let~$\alpha$ be~$2$ if~$N=4$ and a square root of~$2$ if~$N=8$, and we define
\[
O^\prime=D(\alpha^{-1}) \cdot O_n \cdot D(\alpha).
\]
It can now be checked that the valuation of the elements in the~$i^{th}$ column of~$O^\prime$ is at least~$8i/N$.

We will show that the matrix~$O^\prime$ has determinant with valuation~$8/N \cdot n(n+1)$ by showing that it is the product of two matrices, one of which is the diagonal matrix~$D(\alpha)$, and one of which is a matrix with determinant a unit.

We define the matrix~$P$ to be~$D(\alpha)^{-1} \cdot O^\prime$. The entries of~$P$ are given by~$P_{i,j}=\alpha^{-i} \cdot O^\prime_{i,j}$ and are therefore elements of the ring of integers of~$\qtwo (\sqrt{2})$, because the valuation of elements in the~$i^{th}$ column of~$O^\prime$ is at least~$i$. Therefore, we can define the matrix~$P^\prime$ by reducing the entries of~$P$ modulo~$2$; if ~$P^\prime$ has determinant~$1$ in~$\ftwo$ then~$O^\prime$ has determinant~$8/N \cdot n(n+1)$.

By tracing the definition of the columns through all of this, we see that the odd-indexed columns of the mod~$2$ matrix in weight~$(2i+1)\cdot k$ have generating functions
\[
\col(2l+1)=\frac{x^{l+1}y^{2l+1}\cdot (1+x)^i}{(1+x)^{2l+1}},
\]
and the even-indexed columns of the matrix have generating functions
\[
\col(2l)=\frac{x^{l}y^{2l}\cdot (1+x)^i}{(1+x)^{2l} }.
\]

This is clear from the discussion earlier if~$N=4$ or~$N=8$ and~$4$ does not divide the weight. 

If~$4$ does divide the weight, then we consider the quotient~$E^*_{k,\chi}/V^*_{k,\chi}$; notice that the $q$-expansions of the Eisenstein series~$E^*_{2,\chi}$ and~$E^*_{k,\chi}$ are congruent modulo~$2$, by Fermat's little theorem. Let~$z_k:=(E^*_{k,\chi}/V^*_{k,\chi}-1)/\lambda_k$, where~$\lambda_k$ is a square root of~$-2k/B_{k,\chi}$. This is also an overconvergent modular form, so it is an element of~$\qtwo(\sqrt{2})\langle z_8 \rangle$. Then the effect of the changes above on the columns is to send the multiplier~$(1+\lambda_k z_k)$ to~$1+z_k$. We see that~$z_k \equiv z_8 \mod 8$ because
\[
z_k=\frac{E^*_{k,\chi}-V^*_{k,\chi}}{\lambda_k V^*_{k,\chi}} \equiv E^*_{k,\chi}-V^*_{k,\chi} \equiv E^*_{2,\chi}-V^*_{2,\chi} \equiv \frac{E^*_{2,\chi}-V^*_{2,\chi}}{\sqrt{2} V^*_{2,\chi}} \equiv z_8.
\]
Therefore the multiplier~$E^*_{k,\chi}/V^*_{k,\chi}$ when expanded as a rational function of~$z_8$  and then altered by the procedure above is congruent modulo~$2$ to~$1+z_8$.

For example, we consider the~$z_8$ expansion of~$E^*_{4,\chi}/V^*_{4,\chi}$. After sending~$2z_8$ to~$z_8$, we see that we have
\[
\frac{11+z_8+6z_8^2-6z_8^3-z8^4-11z_8^5}{11+6z_8^2-z_8^4} \equiv \frac{1+z_8+z8^4+z8^5}{1+z_8^4} \equiv 1+z_8 \mod 2.
\]

I would like to thank Robin Chapman~\cite{chapman-suggestion} for the idea behind the following proof. To show that the $n \times n$ mod~$2$ matrix has determinant~$1$, we will show that the elements~$C:=\left\{\col^\prime(1)=\col{1} \cdot (1+x)^N,\ldots,\col^\prime(N)=\col{N}\cdot (1+x)^N\right\}$ are a basis 
of the ring~$\ftwo[x]/(x^{N+1})$. Because~$(1+x)$ is a unit in the ring, we may multiply each~$\col(i)$ by~$(1+x)^{N-i}$ to make the calculations easier.
We see that there are exactly the right number of elements, so we must check that they are linearly independent. 

We write~$\sum_{i=1}^N \lambda_i \col^\prime(i)=0$. We will show that all of the~$\lambda_i=0$, so that the columns are linearly independent.

We see that the element~$\col^\prime(1)=x(1+x)^{N-1}$ is the only element of the set~$C$ which has an~$x^N$ term. Therefore~$\lambda_1=0$. Also, we see that~$\col^\prime(2)=x(1+x)^{N-2}$ is now the only nonzero term with an~$x$ term, so therefore~$\lambda_2=0$.

By continuing this process, we show that each~$\lambda_i$ must be zero. Therefore the set~$C$ is composed of linearly independent elements and hence it is a basis. So the determinant of the mod~$2$ matrix is~$1$. Hence we have shown Theorem~\ref{determinant-theorem} and therefore Theorem~\ref{main-theorem}. \qed

\section{Acknowledgements}

I would like to thank Kevin Buzzard for many helpful conversations and much inspiration while he was my PhD advisor. 
I would also like to thank Edray Goins and Ken McMurdy for helpful conversations, and Robin Chapman for his suggestion on a shorter proof of Theorem~\ref{determinant-theorem}. 
I would also like to thank David Burns and Frazer Jarvis for their careful reading of my PhD thesis~\cite{kilford-thesis}, which is where the results of this paper first appeared.

Some of the computer calculations were executed by the computer algebra package {\sc magma}~\cite{magma} running on the machine \textbf{crackpipe}, which was bought by Kevin Buzzard with a grant from the Central Research Fund of the University of London. I would like to thank the CRF for their support. Some other calculations were performed on the machine \textbf{shimura}, which is based at the University of California at Berkeley. I would like to thank the administrators, Wayne Whitney and William Stein, for allowing me to hold an account on their machine.


\end{document}